\documentclass{amsart}
\usepackage{amsmath,amssymb,enumerate,bbm,hyperref}

\newcommand{\assign}{:=}
\newcommand{\emdash}{---}
\newcommand{\nonesep}{}
\newcommand{\tmem}[1]{{\em #1\/}}
\newcommand{\tmop}[1]{\ensuremath{\operatorname{#1}}}
\newcommand{\tmtextit}[1]{{\itshape{#1}}}
\newenvironment{enumeratenumeric}{\begin{enumerate}[1.] }{\end{enumerate}}
\newtheorem{corollary}{Corollary}
\newtheorem{definition}{Definition}
\newtheorem{proposition}{Proposition}
\newtheorem{theorem}{Theorem}

\begin{document}

\title{Constructive pointfree topology eliminates non-constructive
representation theorems from Riesz space theory}\author{Bas
Spitters}\maketitle

\begin{abstract}
  In Riesz space theory it is good practice to avoid representation theorems
  which depend on the axiom of choice. Here we present a general methodology
  to do this using pointfree topology. To illustrate the technique we show
  that almost f-algebras are commutative. The proof is obtained relatively
  straightforward from the proof by Buskes and van Rooij by using the
  pointfree Stone-Yosida representation theorem by Coquand and Spitters.
\end{abstract}

The Stone-Yosida representation theorem for Riesz
spaces~{\cite{RieszI,RieszII}} shows how to embed every Riesz space into the
Riesz space of continuous functions on its spectrum.

\begin{theorem}
  {\tmem{[Stone-Yosida]}} Let $R$ be an Archimedean Riesz space (vector
  lattice) with unit. Let $\Sigma$ be its (compact Hausdorff) space of
  representations. Define the continuous function $\hat{r} (\sigma) \assign
  \sigma (r)$ on $\Sigma$. Then $r \mapsto\hat{r}$ is a Riesz embedding of $R$
  into $C (\Sigma, \mathbbm{R})$.
\end{theorem}

The theorem is very convenient, but sometimes better avoided, since it leads
out of the theory of Riesz spaces. To quote Zaanen~{\cite{Zaanen:Riesz}}:

\begin{quotation}
  Direct proofs, although sometimes a little longer than proofs by means of
  representation [theorems], often reveal more about the situation under
  discussion.
\end{quotation}

Similar concerns where discussed by Buskes, de Pagter and van
Rooij~{\cite{small,Functional-calc}}. They proposed to avoid the use of the
axiom of choice by restricting the size of the Riesz spaces~{\cite{small}}. We
provide a solution based on pointfree topology, which like~{\cite{small}}
avoids also the {\tmem{countable}} axiom of choice, but moreover avoids the
axiom of excluded middle{\footnote{$P \vee \neg P$.}}. This allows the results
to be applied in non-standard contexts. For instance, one can translate a
theorem about one C*-algebra to a theorem about a continuous field of
C*-algebras~{\cite{BM-StoneW,BM:Spectrum,BM:GelfM,BM:Glob}}. In turn, the
results about commutative C*-algebras may be obtained directly using Riesz
spaces~{\cite{CS:Gelfand}}. This is used in applications of topos theory to
quantum theory~{\cite{HLS}}.

Our strategy is as follows. First we replace the topological space of
representations by a locale, a point-free space. This typically removes the
need of the axiom of choice~{\cite{Mulvey:geometry}}. Then we proceed to use
pointfree topology, locale theory using a only basis for the
topology~{\cite{Sambin:1987}}. Since the space of representations of a Riesz
space is compact Hausdorff, it can be described explicitly by the finite
covering relation on the lattice of basic opens. This lattice can be defined
directly using the Riesz space
structure~{\cite{Coquand:Stone,Coquand/Spitters:StoneY}}.

To illustrate the method we will reprove the results by Buskes and van
Rooij~{\cite{BuskesvanRooij}}.

\section{Preliminaries}

\begin{definition}
  A {\tmem{Riesz space}} is a vector space with compatible lattice operations
  {\emdash} i.e.~$f \wedge g + f \vee g = f + g$ and if $f \geqslant 0$ and $a
  \geqslant 0$, then $a f \geqslant 0$. A {\tmem{(strong) unit}} 1 in a Riesz
  space is an element such that for all $f$ there exists a natural number $n$
  such that $| f| \leqslant n \cdot 1$. A Riesz space is {\tmem{Archimedean}}
  if for all $n$, $n |x| \leqslant y$ implies that $x = 0$.
\end{definition}

\subsection{Topologies, locales, lattices}

A topological space may be presented in its familiar set theoretic form, as
such it is a complete distributive lattice of open sets with the operations of
union and intersection. The category of {\tmem{frames}} has as objects
lattices with a finitary meets and infinitary joins such that $\wedge$
distributes over $\bigvee$. Frame maps preserve this structure. A continuous
function $f : X \rightarrow Y$ defines a frame map $f^{- 1} : O (Y)
\rightarrow O (X)$. Since this map goes in the reverse direction it is often
convenient to consider the category of {\tmem{locales}}, the opposite of the
category of frames. In fact, there is a categorical adjunction between the
category of locales and the category of topological spaces. This restricts to
an equivalence of categories for compact Hausdorff spaces and compact regular
locales. In general, the axiom of choice is needed to move from locales to
topological spaces. Hence, by staying on the localic side it is often possible
to avoid the axiom of choice. However, one can go even further. Since compact
locales are determined by the finitary coverings one may restrict one's
attention to the finitary covering relation on a base of the topology. This
base is a normal{\footnote{A lattice is {\tmem{normal}} if for all $b_1, b_2$
such that $b_1 \vee b_2 = \top$ there are $c_1, c_2$ such that $c_1 \wedge c_2
= \bot$ and $c_1 \vee b_1 = \top$ and $c_2 \vee b_2 = \top$. The opens of a
normal topological space form a normal lattice.}} distributive lattice.

For the spectrum of a Riesz space $R$ a base for the topology can be described
very explicitly. Recall that a base for the topology of the spectrum $\Sigma$
is defined by the opens $\{\sigma | \hat{a} (\sigma) > 0\}$. Let $P$ denote
the set of positive elements of $R$. For $a, b$ in $P$ we define $a
\preccurlyeq b$ to mean that there exists $n$ such that $a \leq nb$. The
following proposition is proved in~{\cite{Coquand/Spitters:StoneY}} and
involves only elementary considerations on Riesz spaces.

\begin{proposition}
  \label{nine}$L (R) : = (P, \vee, \wedge, 1, 0, \preccurlyeq)$ is a
  distributive lattice. In fact, if we define $D : R \rightarrow L (R)$ by $D
  (a) \assign a^+$, then $L (R)$ is the free lattice generated by $\{D (a) | a
  \in R\}$ subject to the following relations:
  \begin{enumeratenumeric}
    \item $D (a) = 0$, if $a \leq 0$;

    \item $D (1) = 1$;

    \item $D (a) \wedge D (- a) = 0$;

    \item $D (a + b) \leqslant D (a) \vee D (b)$;

    \item $D (a \vee b) = D (a) \vee D (b)$.
  \end{enumeratenumeric}
\end{proposition}

We have $D (a) \leqslant D (b)$ if and only if $a^+ \preccurlyeq b^+$ and $D
(a) = 0$ if and only if $a \leqslant 0$. We write $a \in (p, q) \assign (a -
p) \wedge (q - a)$ and observe that this is an element of $R$.

For $a$ in $R$ we define its norm $\|a\|= \inf \{q | a \leqslant q 1\}$.

The corresponding locale{\footnote{From this point onwards $\Sigma$ is the
spectrum considered as a locale. If we want to treat it as a topological space
we write $\tmop{pt} \Sigma$.}} (complete distributive lattice) $\Sigma$ is the
one defined by the same generators and relations together with the relation $D
(a) = \bigvee_{s > 0} D (a - s)$. The generators and relations above may also
be read as a propositional geometric theory~{\cite{Vic:LocTopSp}} by reading
$\leqslant$ as $\Rightarrow$. A model $m$ of this theory defines a
representation $\sigma_m$ of the Riesz space by
\[ \sigma_m (a) \assign (\{q|m \models D (q \cdot 1 - a)\}, \{q|m \models D (a
   - q \cdot 1)\}), \]
where the right hand side is a Dedekind cut in the rationals and hence a real
number. Such a $\sigma_m$ is a point of the locale $\Sigma$. This motivates
the interpretation of $D (a)$ as $\{\sigma | \hat{a} (\sigma) > 0\}$: the
models which make the proposition $D (a)$ true coincide with the points
$\sigma$ such that $\hat{a} (\sigma) > 0$. Proving that there are enough such
models/points requires the axiom of choice. We avoid this axiom by working
with the propositions/opens instead.

\begin{theorem}
  \label{Loc:StY}{\tmem{[Localic Stone-Yosida]}} The map $\hat{\cdot} : R
  \rightarrow \tmop{Loc} (\Sigma, \mathbbm{R})$ defined by the frame map
  $\hat{a} (p, q) \assign a \in (p, q)$ is a norm-preserving Riesz morphism.
  Its image is dense with respect to the uniform topology on $\tmop{Loc}
  (\Sigma, \mathbbm{R})$.
\end{theorem}

\begin{proof}
  The map $\hat{\cdot} : R \rightarrow \tmop{Loc} (\Sigma, \mathbbm{R})$ is
  norm-preserving; see~{\cite{Coquand/Spitters:StoneY}}. It remains to prove
  the density. For this consider a natural number $N$ and a continuous $f$ on
  $\Sigma$ such that $0 \leqslant f \leqslant 1$. We need to find an element
  $a$ of $R$ such that $\hat{a}$ is close to $f$. The set $\bigcup^N_{k = 0} f
  \in ((k - 1) / N, (k + 1) / N)$ covers $\Sigma$. By Proposition~3.1
  of~{\cite{Coquand:Stone}} there exists a partition of unity $p_i$ in the
  Riesz space such that $\sum p_i = 1$ and $D (p_i)$ is contained in some open
  $f \in ((k_i - 1) / N, (k_i + 1) / N)$ in $\Sigma$. Concretely, $p_i
  \preccurlyeq (f - (k_i - 1) / N) \wedge ((k_i + 1) / N - f)$. Consequently,
  \[ |f - \sum k_i \widehat{p_i} | = |f \sum \widehat{p_i} - \sum k_i
     \widehat{p_i} | = | \sum (f - k_i) \widehat{p_i} | \leqslant \frac{1}{N}
     . \]
\end{proof}

The map $\hat{\cdot}$ is a Riesz embedding if $R$ is Archimedean.

\begin{corollary}
  There is a norm-preserving Riesz morphism of $R$ into an f-algebra such that
  the image is dense.
\end{corollary}

The axiom of choice implies that compact regular locales have enough points
and hence we obtain the more familiar formulation of the theorem by working
with the topological space $\tmop{pt} \Sigma$ of the points of the spectrum.
However, in practice, only the localic version is needed.

\begin{corollary}
  {\tmem{[Stone-Yosida]}} The map $\hat{\cdot} : R \rightarrow C (\tmop{pt}
  \Sigma, \mathbbm{R})$ defined by the frame map $\hat{a} (p, q) \assign a \in
  (p, q)$ is a norm-preserving Riesz morphism. Its image is dense with respect
  to the uniform topology on $C (\tmop{pt} \Sigma, \mathbbm{R})$.
\end{corollary}

\section{The results}

\begin{definition}
  An {\tmem{almost $f$-algebra}} is a Riesz space with multiplication such
  that $a \cdot b \geqslant 0$ if $a, b \geqslant 0$, and $a \wedge b = 0$
  implies $a \cdot b = 0$.
\end{definition}

If $E$ is a Riesz space, a bilinear map A of $E \times E$ into a vector space
F is called {\tmem{orthosymmetric}} if
\[ f \wedge g = 0 \Rightarrow \nonesep A (f, g) = 0 \]
for all $f, g \in E$.

A \tmtextit{partition of unity} is a list $u_i$ such that $\sum u_i = 1$ and
$0 \leqslant u_i \leqslant 1$. If $u, v$ are partitions of unity in an almost
f-algebra, then so is $u \cdot v$: $\sum_i \sum_j u_i v_j = 1 \cdot \sum_j
v_j$.

\begin{theorem}
  \label{Thm:main}Let E be Riesz spaces with unit and let $F$ be Archimedean
  and let A be a orthosymmetric positive bilinear map $E \times E \rightarrow
  F$. Let $\bar{E}$ be an f-algebra in which $E$ is dense and let $F'$ be the
  uniform completion of $F$. Then $A$ extends uniquely to a orthosymmetric
  positive bilinear map from $\bar{E} \times \bar{E}$ to $F'$ and $A (f, g) =
  A (1, f g)$ for all $f, g$ in $E$.
\end{theorem}

\begin{proof}
  Let $f, g$ be in $E$.
  \[ A (f, g) = A (f^+, g^+) + A (f^-, g^-) - A (f^-, g^+) - A (f^+, g^-) . \]
  So, it suffices to consider the case $0 \leqslant f, g \leqslant 1$. Let $k$
  be a natural number. Define $u_n \assign k (f \vee \frac{n}{k} \wedge
  \frac{n + 1}{k})$, whenever $0 \leq n < k$. Define $v_0 \assign 1 - u_0$ and
  $v_n \assign u_n - u_{n + 1}$ and $v_k \assign u_k .$ The set $\{v_0,
  \ldots, v_k \}$ is a partition of unity {\emdash} that is, $\sum v_i = 1$
  and $0 \leqslant v_i \leqslant 1$. Moreover, $v_n \perp v_m$, whenever $|n -
  m| > 1$ and such that $|fv_n - \frac{n}{k} v_n | \leq \frac{1}{k}$. By
  repeating a similar construction for $g$ we find a partition of unity $v'$.
  Then $w_{i j} \assign v_i v'_j$ is again a partition of unity. For
  convenience, we reindex $w$ by one natural number to obtain a sequence
  $w_n$. We define $\alpha_n, \beta_n$ such that $|f - \sum_n \alpha_n w_n | <
  \frac{1}{k}$ and $|g - \sum_n \beta_n w_n | < \frac{1}{k}$.

  Let $\varepsilon = \frac{1}{k}$. Set $f' : = \sum \alpha_n w_n$, $g' : =
  \sum \beta_n w_n$ and $h' : = \sum \alpha_n \beta_n w_n$. Then
  \begin{eqnarray}
    |A (f, g) - A (f', g') | & = & |A (f - f', g) + A (f', g - g') |
    \nonumber\\
    & \leqslant & \varepsilon A (1, 1) + \varepsilon A (1, 1) \nonumber
  \end{eqnarray}
  since $|f - f' |, |g - g' | \leqslant \varepsilon$ and $A$ is positive.
  Thus, it suffices to show that
  \[ |A (f', g') - A (1, h') | \leqslant 2 \varepsilon A (1, 1) . \]
  Observe that for all $n, m$ in $\{1, . . ., N\}$,

  \ \ \ \ \ if $|n - m| > 1$, then $w_n \perp w_m$, so $A (w_n, w_m) = 0$;

  \ \ \ \ \ if $|n - m| \leqslant 1$, then $| \alpha_n - \alpha_m | \leqslant
  2 \varepsilon$.

  It follows that
  \begin{eqnarray}
    |A (f', g') - A (1, h') | & = & | \sum \alpha_n \beta_m A (w_n, w_m) -
    \sum_{n, m} \alpha_m \beta_m A (w_n, w_m) | \nonumber\\
    & \leqslant & \sum | \alpha_n - \alpha_m | | \beta_m |A (w_n, w_m)
    \nonumber\\
    & \leqslant & 2 \varepsilon \sum_{n, m} A (w_n, w_m) = 2 \varepsilon A
    (1, 1) \nonumber
  \end{eqnarray}
  The last inequality follows from the observation above and the inequality $|
  \beta_m | \leqslant 1$.

  Changing the roles of the $\alpha$s and $\beta$s we have that $|A (g', f') -
  A (1, h') | \leqslant 2 \varepsilon A (1, 1)$. Hence $|A (f', g') - A (g',
  f') | \leqslant 4 \varepsilon A (1, 1)$ and $|A (f, g) - A (g, f) |
  \leqslant 8 \varepsilon A (1, 1)$. Finally, $|h' - f' g' | \leqslant | \sum
  \alpha_m \beta_m w_n w_m - \sum \alpha_n \beta_m w_n w_m | \leqslant \sum |
  \alpha_n - \alpha_m | | \beta_m |w_n w_m \leqslant 2 \varepsilon$. Hence
  $|h' - f g|$ is small. Since $F$ is Archimedean, $A (f, g) = A (1, f g)$.
\end{proof}

The completion of $E$ and $F$ in the previous proof are used to define the
multiplication on $E$ and to be able to extend $A$ to this completion of $E$.
This, however, can be avoided as follows. The joint partition of unity can be
obtained as in Theorem~\ref{Loc:StY}. The proof above then shows that for each
$\varepsilon$, $|A (f, g) - A (g, f) | \leqslant 8 \varepsilon A (1, 1)$. This
implies the following result.

\begin{corollary}
  \label{Cor:comm}Let E and F be Riesz spaces of which $F$ is Archimedean.
  Let A be an orthosymmetric positive bilinear map $E \times E \rightarrow F$.
  Then
  \[ A (f, g) = A (g, f) \hspace{2em} (f, g \in E) . \]
\end{corollary}

\begin{proof}
  Take $f, g \in E$. Let $E_0$ be the Riesz subspace of $E$ generated by $\{f,
  g\}$. Then $|f | + |g |$ is a unit in $E_0$. Without restriction, suppose
  that $E_0$ is $E$. The result is now follows from the previous theorem.
\end{proof}

We have proved the following result in an entirely explicit way by a
straightforward analysis of the proof by Buskes and van Rooij.

\begin{corollary}
  Every Archimedean almost f -algebra is commutative.
\end{corollary}

\section{Internal real numbers}

Buskes and van Rooij use the Stone-Yosida representation theorem combined with
Dini's theorem to show that a certain sequence of elements in a Riesz space
converges. This can be replaced by applying the following result which does
not require the sequence to be decreasing.

\begin{theorem}
  \label{sheaf}Let $e_n$ be a sequence of expressions in the language of Riesz
  spaces such that $e_n$ converges {\tmem{constructively}} when interpreted in
  the Riesz space of real numbers. Then $e_n$ converges uniformly when
  interpreted in any Riesz space with strong unit.
\end{theorem}

\begin{proof}
  The Riesz space can be (densely) embedded into a space $C (\Sigma)$ and
  hence its elements may be interpreted as global sections of the real number
  object in the topos $\tmop{Sh} (\Sigma)$ of sheaves over
  $\Sigma$~{\cite{johnstone02b,maclanemoerdijk92}}. Now, if $a_n$ converges to
  0 in the internal language of $\tmop{Sh} (\Sigma)$. Then for each $q$ there
  exists $n$ such that $a_n \leqslant q$ internally. This is interpreted as:
  for each $q$ there exists a (finite) cover $U_i$ of $\Sigma$ and $n_i$ such
  that $a_{n_i} \leqslant q$ on $U_i$. Taking $n = \min n_i$ we see that $a_n
  \leqslant q$ on $\Sigma$.
\end{proof}

Sheaf theory may seem to be a very complex tool to use for such a simple
lemma, however, when applied in concrete cases we obtain natural results. For
instance, Buskes and van Rooij apply Dini's theorem and the Stone-Yosida
representation theorem to prove that the sequence $[(f \wedge g) h - n f (g
\wedge h)]^+$ converges uniformly. We first prove this for the Riesz space of
real numbers. Fix $m$ in $\mathbbm{N}$. We may assume that $f, g, h \leqslant
1$. Moreover, either{\footnote{The case distinction $f \geqslant \frac{1}{m}$
or $f \leqslant \frac{1}{m}$ is not constructive/continuous.}} $f \geqslant
\frac{1}{m}$ or $f \leqslant \frac{2}{m}$. We may assume that $f \geqslant
\frac{1}{m}$ and similarly that $g, h \geqslant \frac{1}{m}$. Choosing $n =
m^2$ shows that $(f \wedge g) h \leqslant 1$ and $n f (g \wedge h) \geqslant n
\cdot \frac{1}{m} \cdot \frac{1}{m} \geqslant 1$. Hence if $n \geqslant m^2$,
then $[(f \wedge g) h - n f (g \wedge h)]^+ \leqslant \frac{2}{m}$ in all
cases. The interpretation of this statement in the sheaf model $\tmop{Sh}
(\Sigma)$ defined from the spectrum $\Sigma$ of a Riesz space is: there is a
finite cover $U_i$ of $\Sigma$ such that $[(f \wedge g) h - n f (g \wedge
h)]^+ \leqslant \frac{2}{m}$ is true on each $U_i$. A finite cover gives rise
to a partition $u_i$ of unity such that $D (u_i) \subset U_i$. So, that $u_i
[(f \wedge g) h - n f (g \wedge h)]^+ \leqslant \frac{2}{m}$ and hence
\[ [(f \wedge g) h - n f (g \wedge h)]^+ = \sum u_i [(f \wedge g) h - n f (g
   \wedge h)]^+ \leqslant \frac{2}{m} . \]

Takeuti's use of Boolean valued models~{\cite{Takeuti}} to obtain non-standard
results from familiar theorems has a similar flavor as Theorem~\ref{sheaf}.
Boolean valued models are a special class of sheaf models.

\section{Conclusion}

We have illustrated how the use of locale theory, presented by a normal
distributive lattice of basic elements, naturally translates proofs which
depend on the axiom of choice to simpler lattice theoretic proofs which avoid
the axiom of choice, even in its countable form, and the principle of excluded
middle. Buskes and van Rooij had previously proposed different methods to
avoid the axiom of choice. An advantage of our approach is that it is valid in
any topos. It also provides a logical tool to remove the use of representation
theorems from Riesz space theory, the importance of avoiding representation
theorems was stressed by Zaanen.

\section{Acknowledgments}

I would like to thank Gerard Buskes, Thierry Coquand and Arnoud van Rooij for
comments on an earlier draft of this paper.

\end{document}